# A study of divergence from randomness in the distribution of prime numbers within the arithmetic progressions $1 + 6n$ and $5 + 6n$

Andrea Berdondini

ABSTRACT. In this article I present results from a statistical study of prime numbers that shows a behaviour that is not compatible with the thesis that they are distributed randomly. The analysis is based on studying two arithmetical progressions defined by the following polynomials: $(1 + 6n, 5 + 6n, n \in \mathbf{N})$ whose respective numerical sequences have the characteristic of containing all the prime numbers except 3 and 2. If prime numbers were distributed randomly, we would expect the two polynomials to generate the same number of primes. Instead, as the reported findings show, we note that the polynomial $5 + 6n$ tends to generate many more primes, and that this divergence grows progressively as more prime numbers are considered. A possible explanation for this phenomenon can be found by calculating the number of products that generate composite numbers which are expressible by the two polynomials. This analysis reveals that the number of products that generate composite numbers expressible by the polynomial $1 + 6n$ is $(n + 1)^2$, while the number of products that generate composites expressible by the polynomial $5 + 6n$ is $n(n + 1)$, with $n \in \mathbf{N}$. As a composite number is a non-prime number, this difference incited me to analyse the distribution of prime numbers generated by the two polynomials. The results, based on studying the first (approx.) 500 million prime numbers, confirm the fact that the number of primes that can be written using the polynomial $1 + 6n$ is lower than the number of primes that can written using the polynomial $5 + 6n$, and that this divergence grows progressively with the number of primes considered.

**Introduction**

The study of the distribution of prime numbers has always occupied a fundamental place in mathematics and for this reason any experimental evidence of their divergence from randomness represents an important finding.

In a recent work, Robert Lemke Oliver and Kannan Soundararajan [1], reported statistical data on an anomaly in the distribution of pairs comprising two consecutive primes. Another important piece of statistical data regarding prime numbers has been found by studying the distribution of increments (difference between two consecutive distances) Pradeep Kumar, Plamen Ch. Ivanov, H. Eugene Stanley [2], which was found to have statistical characteristics that would be highly improbable if the distribution of prime numbers were random.

In this analysis, I use three arithmetical progressions with a common difference of 6, defined by the three following polynomials: $1 + 6n, 3 + 6n, 5 + 6n$ ($n \in \mathbf{N}$), to generate all odd numbers. Of these three polynomials, my analysis will concentrate on just two: $1 + 6n, 5 + 6n$ which have the characteristic of generating all the prime numbers with the exceptions of 2 and 3. In fact, the discarded polynomial $3 + 6n$ does no more than generate all odd numbers divisible by 3. In the following section I will demonstrate a very important finding which sparked my curiosity, leading me to carry out a statistical analysis of the difference in the number of primes that can be generated using the two polynomials.

The finding in question concerns the calculation of the number of products whose results are composite numbers that can be expressed using the two polynomials under consideration. From this analysis, as I demonstrate in the following section, it can be seen that $(n + 1)^2$ products have as their result composite numbers expressible by the polynomial $1 + 6n$, while $n(n + 1)$ products have as their result composite numbers expressible by the polynomial $5 + 6n$, with $n \in \mathbf{N}$.

As a composite number is a non-prime number, this difference led me to think that prime numbers are more numerous within the group of odd numbers produced by the polynomial $5 + 6n$.

This finding proves to be very significant because it contradicts the thesis that the distribution of prime numbers is random. If this were indeed the case, then prime numbers should not show any preference between the two polynomials.

**Theory**

Let us begin by defining three arithmetical progressions with a common difference of 6 which generate, within natural numbers, the whole subset of odd numbers. These arithmetical progressions are defined by the following polynomials: $1 + 6n$, $3 + 6n$, $5 + 6n$, where $n \in \mathbf{N}$. We thereby obtain the whole sequence of odd numbers.

|           | $1 + 6n$ | $3 + 6n$ | $5 + 6n$ |
|-----------|----------|----------|----------|
| $n = 0$   | 1        | 3        | 5        |
| $n = 1$   | 7        | 9        | 11       |
| $n = 2$   | 13       | 15       | 17       |
| $n = \cdots$ | ...   | ...      | ...      |

Let us now consider just the polynomial $3 + 6n$, which can also be written in this way: $3(1 + 2n)$; one can see straight away that this polynomial generates just one prime number, 3. Accordingly, the polynomials $1 + 6n$ and $5 + 6n$ generate all the prime numbers with the exceptions of 2 and 3.

Let us now analyse all the products whose results are composite numbers that can be expressed using the two polynomials $1 + 6n$ and $5 + 6n$. Possible cases are represented by three combinations, the first two of which represent the multiplication of numbers within the same arithmetic progression. In the third case, the multiplication is of two numbers, one of which belongs to the arithmetic progression defined by the polynomial $1 + 6n$ and other belonging to the arithmetic progression defined by the polynomial $5 + 6n$. These products define all the composite numbers present within the two arithmetical progressions considered.

The three possible products are:

$$(1 + 6n_1)(1 + 6n_2) = 1 + 6n_1 + 6n_2 + 36n_1n_2$$
$$(5 + 6n_1)(5 + 6n_2) = 25 + 30n_1 + 30n_2 + 36n_1n_2$$
$$(1 + 6n_1)(5 + 6n_2) = 5 + 30n_1 + 6n_2 + 36n_1n_2$$

I shall now demonstrate that the first two equations can be written using the polynomial $1 + 6n$ and the third using the polynomial $5 + 6n$.

$$1 + 6n_1 + 6n_2 + 36n_1n_2 = 1 + 6(n_1 + n_2 + 6n_1n_2)$$

Defining $n_3 = n_1 + n_2 + 6n_1n_2$ we have:

$$1 + 6(n_1 + n_2 + 6n_1 n_2) = 1 + 6n_3$$

$$25 + 30n_1 + 30n_2 + 36n_1 n_2 = 1 + 6(4 + 5n_1 + 5n_2 + 6n_1 n_2)$$

Defining $n_3 = 4 + 5n_1 + 5n_2 + 6n_1 n_2$ we obtain:

$$1 + 6(4 + 5n_1 + 5n_2 + 6n_1 n_2) = 1 + 6n_3$$

$$5 + 30n_1 + 6n_2 + 36n_1 n_2 = 5 + 6(5n_1 + n_2 + 6n_1 n_2)$$

Defining $n_3 = 5n_1 + n_2 + 6n_1 n_2$ we obtain:

$$5 + 6(5n_1 + n_2 + 6n_1 n_2) = 5 + 6n_3$$

From which we may deduce the following rule: when two odd numbers that have both been generated by the same polynomial ($1 + 6n$ or $5 + 6n$) are multiplied, one obtains a number that can always be expressed by the polynomial $1 + 6n$. On the other hand, when two odd numbers, one of which has been generated by the polynomial $1 + 6n$ and the other by the polynomial $5 + 6n$ are multiplied, one obtains an odd number that can be expressed by the polynomial $5 + 6n$.

Using this simple rule, it is an easy matter to calculate the number of products that are generated by multiplying together two odd numbers which have in turn been generated by the two polynomials under consideration. Thereby, the number 1 is excluded. This is because, as we want to analyse the distribution of prime numbers, we want to consider only the products that generate non-prime odd numbers.

Let us start from the products whose results can be written using the polynomial $1 + 6n$. As demonstrated above, these represent the products of numbers generated by the same polynomial. As we have two polynomials ($1 + 6n$, $5 + 6n$), we have two sequences of odd numbers:

|  | $1 + 6n$ | $5 + 6n$ |
|---|---|---|
| $n = 0$ | 1 | 5 |
| $n = 1$ | 7 | 11 |
| $n = 2$ | 13 | 17 |
| $n = \cdots$ | ... | ... |

Let us now turn to the polynomial $1 + 6n$ and let us calculate all the possible products of two numbers generated by it for a fixed $n$, therefore taking into consideration the numbers of the numerical sequence up as far as the value $1 + 6n$. Let us begin from the first term (1, with $n = 0$) that can be multiplied by itself and by all the other numbers of the sequence up to the term $1 + 6n$. The first term 1 thus produces $n + 1$ products.

Let us move on to the second term: (7, with $n = 1$) that can be multiplied by itself and by all the other numbers of the sequence greater than itself up to the term $1 + 6n$. The second term therefore generates $n$ products. Let us repeat this procedure for each number in the sequence until reaching the term $1 + 6n$. From which we derive, having defined a value for $n$, that the total number of products is:

$$\sum_{p=0}^{n+1} p$$

Removing from the calculation products between the number 1 and other numbers of the arithmetic progression, we obtain:

$$\sum_{p=0}^{n} p$$

Now let us consider the sequence of numbers generated by the polynomial $5 + 6n$. In this case, as the number 1 does not occur in the sequence of odd numbers generated, the number of all possible products between two numbers generated by it for a fixed $n$ is equal to:

$$\sum_{p=0}^{n+1} p$$

By adding the two equations together, we obtain the number of products whose results can be expressed by the polynomial $1 + 6n$. Having set an $n$, one obtains:

$$\sum_{p=0}^{n} p + \sum_{p=0}^{n+1} p$$

Knowing that the sum of the first $n$ odd numbers is equal to $n^2$ we have:

$$\sum_{p=0}^{n} p + \sum_{p=0}^{n+1} p = (n+1)^2$$

Finally, we calculate the number of products whose results can be expressed by the polynomial $5 + 6n$. Using the rule demonstrated above, the polynomial $5 + 6n$ generates the products of two numbers, one of which is generated by the polynomial $1 + 6n$ and the other by the polynomial $5 + 6n$. You can appreciate that in this case the calculation of the number of these products is trivial. In fact, starting from the two sequences generated by the two polynomials:

|  | $1 + 6n$ | $5 + 6n$ |
|---|---|---|
| $n = 0$ | 1 | 5 |
| $n = 1$ | 7 | 11 |
| $n = 2$ | 13 | 17 |
| $n = \cdots$ | ... | ... |

every number generated by the polynomial $5 + 6n$ will now be multiplied by every number generated by the polynomial $1 + 6n$ (excepting the first term, which is 1) in order to obtain all the other remaining composite numbers. Having set an $n$, every number generated by $5 + 6n$ can, then, determine $n$ products. As the numbers generated by the polynomial $5 + 6n$ are equal to $n + 1$, the total number of possible products is:

$n(n + 1)$

As the product of any two odd numbers apart from 1 is a non-prime number, one may assume from the calculation of combinations shown above that the distribution of prime numbers is not uniform across the two polynomials. In the next section I shall present statistical data in which the discrepancy is calculated between the numbers of primes that can be expressed by the two polynomials. These results will clearly show that this difference grows as we increase the number of primes that we take into consideration.

**Results**

I shall now proceed to demonstrate the findings concerning the difference in the numbers of primes within the two arithmetical progressions defined by the two polynomials $1 + 6n$, $5 + 6n$, with $n \in \mathbf{N}$. I shall, then, define a function I shall call $\Delta P(Np)$, which will be equal to the difference between the numbers of primes that can be written with the polynomial $5 + 6n$ and the number of primes that can be written with the polynomial $1 + 6n$, when considering $Np$ prime numbers.

Figure 1 shows the function $\Delta(Np)$, in which the first 489,736,000 prime numbers have been analysed. It can be seen that, even though there are wide variations in the value of $\Delta(Np)$, the difference in the numbers of primes generated by the two polynomials under consideration tends to grow as the number of primes analysed increases.

Figure 2 shows $\Delta P(Np)$ for a lower number of analysed primes, ($Np = 50,000$); in this way one can see more clearly the way it develops in values close to zero. It should be noted that the function $\Delta(Np)$ does not assume negative values.

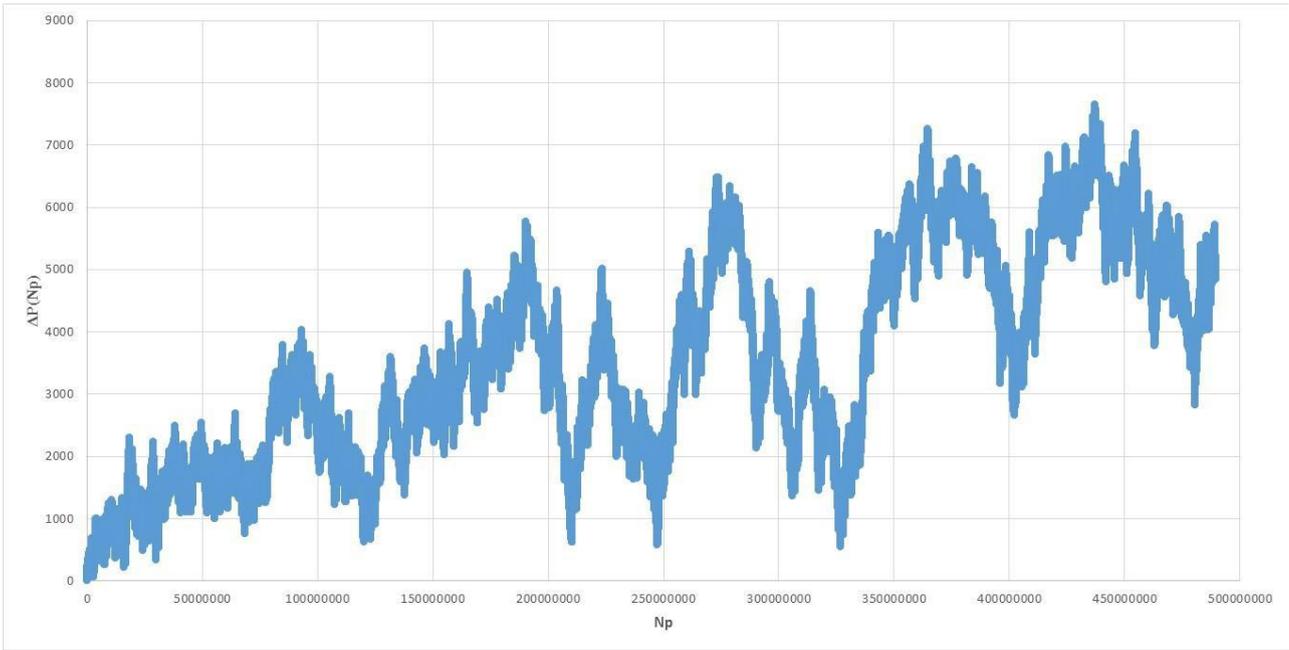

**FIG. 1**: The function $\Delta(Np)$ in the case where the first 489,736,000 prime numbers have been analysed.

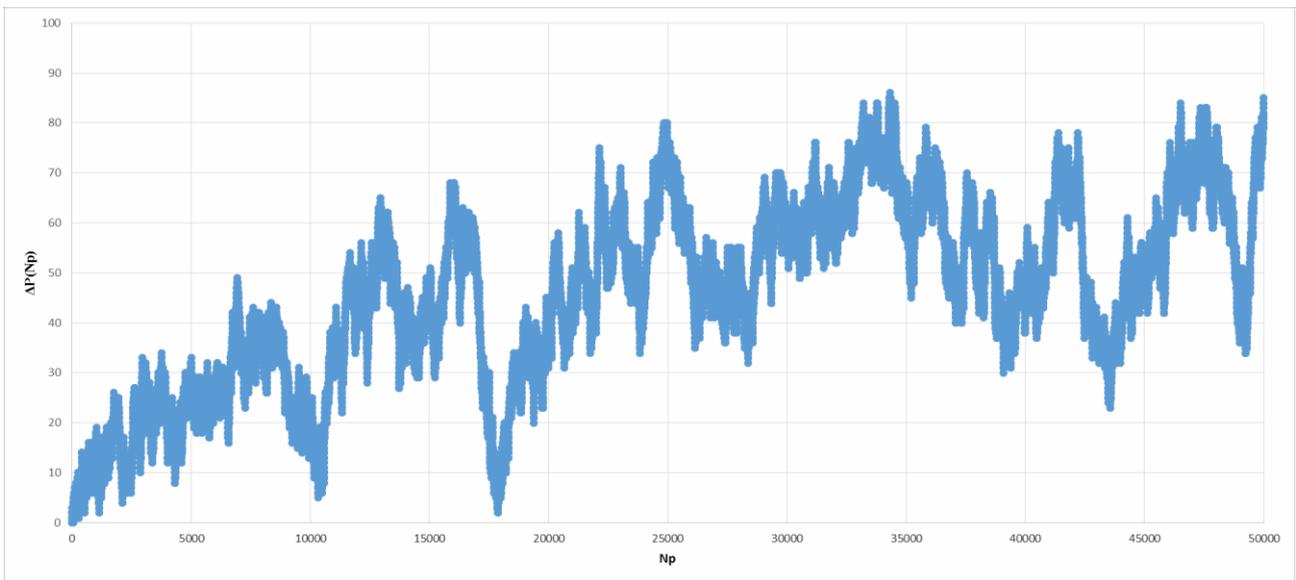

**FIG. 2**: The function $\Delta(Np)$ in the case where the first 50,000 prime numbers have been analysed.

## Conclusion

In this article I have reported data from a statistical analysis that highlights a divergence from randomness in the distribution of prime numbers. This statistical anomaly is manifested by a comparative analysis of the distribution of prime numbers within two degree-one polynomials ($1 + 6n$, $5 + 6n$). These two polynomials define two arithmetical progressions with a common difference of 6, which have the characteristic of generating all the prime numbers with the exceptions of 3 and 2. If one assumes the thesis that prime numbers are distributed randomly, one would expect prime numbers to show no preference between being generated by one polynomial rather than another. The statistical data reported, based on an analysis of the first 500 million prime numbers approximately, shows how the polynomial $1 + 6n$ generates fewer prime numbers than the polynomial $5 + 6n$ and how this difference increases progressively with the amount of data analysed. One possible explanation can be found by calculating the number of products, whose results are composite numbers, that can be expressed by the two polynomials. From this analysis it is demonstrated that the number of products that generate composite numbers expressible by the polynomial $1 + 6n$ is $(n + 1)^2$, while the number of products that generate composites expressible by the polynomial $5 + 6n$ is $n(n + 1)$, with $n \in \mathbf{N}$. As a composite number is a non-prime number, this difference may be considered to be of significance in interpreting the data. Unfortunately, as the values of the composite numbers obtainable by multiplication may involve more than one pair of numbers belonging to both of the arithmetical progressions under consideration, this combinatorial finding cannot be taken as a demonstration in explaining the statistical divergence encountered. For this reason, one can define only one possible conjecture about the fact that the value of the function $\Delta(Np)$ tends to infinity as the number of primes under consideration $Np$ tends to infinity.

One possible interpretation of this statistical finding may be to see it as the formation of an intrinsic ordering between the numbers when they interact with each other through multiplication. This possible divergence from randomness in the distribution of composite numbers consequently involves an ordering among prime numbers.

*E-mail address*: andrea.berdondini@libero.it